**The Purpose of Mathematics according to Plato and Augustine**
**January 9, 2024**

**Douglas J. Dailey**
**Associate Professor of Mathematics**
**Christendom College**
**Front Royal, VA**
**Douglas.Dailey@christendom.edu**

## 1. Introduction

The history of mathematics betrays a sort of inevitability in the development of its theories. There are many examples of cases where one mathematician made a discovery only to find that the same exact discovery had been made by another mathematician working independently. Thus, it was that calculus was founded as a field in similar manners by Newton and Leibniz around the same time. Similarly, Lobachevsky and Bolyai independently founded Non-Euclidean geometry in the early nineteenth century. Only later did Gauss say that he had laid the foundations himself around ten years before either of them, though he never published his results.

The fact of inevitable mathematical discoveries is the starting point for a lecture by a leading Russian mathematician by the name of Igor Rostislavovich Shafarevitch. As a professor at Moscow State University, Shafarevitch was a leading algebraic geometer of the twentieth century. Shafarevitch also stands out politically in Russia as a long-standing dissident of the Soviet government. In the 1970's, he was a member of a group of dissidents which also included Aleksandr Solzhenitsyn. As a result of this activity, he lost his university posting in 1975.[1] His lecture is entitled 'On Certain Tendencies in the Development of Mathematics,' and he delivered it to the Göttingen Academy of Sciences in 1973.

Shafarevitch likens the experience of independent development of a mathematical theory to the players in an orchestra where the theme moves from one instrument to another seamlessly. To him, in order for mathematics to develop in this manner, there must be some purpose or aim that drives mathematical discoveries. Indeed, 'the internal logic of its development much more resembles the work of a single intellect developing its thought in a continuous and systematic way' according to some purpose.[2] Assuming that this is the case, what can we say that the purpose of mathematics is?

In physics, beginning with Newton, the aim was to have a system in which all physical phenomena from the large-scale to the minute were contained and explained by a few simple laws. The complexity of the physical universe and later discoveries such as the development of the theories of electromagnetism and quantum physics have shown that such a system is difficult to determine. Many physicists doubt that such a system will ever be found if, in fact, they have even paused to consider whether such a thing should be pursued.

---

[1] Kishkovsky, Sophia. "Igor Shafarevich, Russian Mathematician With a Mixed Political Legacy, Dies at 93." *New York Times,* March 13, 2017. https://www.nytimes.com/2017/03/13/world/europe/igor-shafarevich-dead-dissident-mathematician.html

[2] Shafarevitch, I.R. 'On Certain Tendencies in the Development of Mathematics.' *Poetics Today* 3, No. 1 (Winter 1982), p. 5.

Unrealistic as the goal of a single unifying system for physics may be, mathematicians, according to Shafarevitch, have never formulated a single aim or goal in the way that physicists did. Without a goal, the only 'ideal it is left with is uncontrolled growth…in all directions'.[3] One can certainly see this proliferation of results in the literature. The online pre-print repository arXiv.org shows that there were over 45,000 mathematics papers posted on their site in 2022 alone. The danger that Shafarevitch points out is that with such proliferation one risks losing all meaning behind the endeavor, and this is especially dangerous in mathematics, a discipline in which order and beauty clearly shine forth in the profound examples of classical and modern results. Mathematics without an aim is akin to the 'inherent contradiction in the concept of a symphony which goes on forever'.[4] Expanding his scope, Shafarevitch observes that all human activity without an aim becomes meaningless.

So what is a mathematician to do? Experience dictates that the position of having no aim is untenable, so people find an aim by way of substitution, or, more precisely, people borrow aims from other areas. Shafarevitch posits two such directions in which mathematicians may borrow an aim. The first is through the practical applications of mathematics. While it is true that practical utility of mathematics is vast and deep, it is 'not the case that these applications inspired its most beautiful achievements'.[5] Rather, we see that logical proof is the means by which the discipline develops and the impetus behind the early accomplishments of mathematics. But neither can we say that logic or rigor is the aim of mathematics, for the goal of mathematics cannot be subordinate to its method.

For Shafarevitch, then, mathematicians need to borrow their aim from a higher sphere, and he says that the only aim that will do is that of religion. To make his point, Shafarevitch goes back to the earliest days of mathematics in which the Pythagoreans are purported to have followed exactly this path in mathematics.

To modern eyes, such a claim may seem shocking. Mathematicians infrequently delve into such matters as the philosophy of mathematics, and even if they do, they return quickly to the more scientific aspect of their discipline. Even for a Christian mathematician, the suggestion that mathematics derives its aim as a religious undertaking may seem like a stretch or, if not, at least a thing difficult to articulate. So what are we to make of this claim?

The outline that Shafarevitch follows in his lecture is strongly influenced by early Greek philosophers, particularly those within the Platonic tradition. Thus, to make sense of Shafarevitch's claim, we would do well to understand the purpose that mathematics served within a person's intellectual formation according to Plato. To place Plato's view into a Christian perspective, we will then investigate the thought of St. Augustine of Hippo, the great fifth century theologian and bishop. Augustine's insight on the role that number plays in the development of reason sheds light on how knowledge of mathematics conduces to knowledge of God.

---

[3] Ibid, p. 7.
[4] Ibid.
[5] Ibid, p. 8.

## 2. Plato

In Plato's *Republic,* Socrates sets out his vision for the ideal city as well as his conception of the education and preparation of its future rulers. A large component of their education is designated for mathematical studies. The specific nature of the mathematical subjects themselves and the place they occupy in the rulers' overall education give an indication for how Socrates views the field of mathematics within the hierarchy of human knowledge. In each example, we will see how utility and even logic are mere 'indirect effects' of mathematical learning as opposed to their main purpose, which is leading the soul towards being and truth.

In her 2020 Aquinas Lecture delivered at Marquette University, Sarah Broadie states that it is uncontroversial among scholars that Plato holds that 'having the right attitude towards truth needed for being…a good ruler presupposes a deep mathematical education'[6]. This stance is evinced in the fact that for the future rulers of Socrates' ideal city, mathematical studies are to be a lifelong companion. Even from a young age 'calculation and geometry…should be presented to the mind in childhood' albeit under no compulsion, but rather as a 'sort of amusement' in order to find out the students' 'natural bent'.[7] The mathematical studies thus begun in an informal manner will later be presented formally to worthy students selected at the age of twenty when they will learn 'to see the natural relationship of [the topics] to one another and to true being'.[8] After ten years of study of the mathematical sciences, students will be prepared to begin studying dialectic, the entire end to which their education is directed.[9]

What is the nature of dialectic, and how does mathematics fit into it? This question is a topic of much scholarly debate, and the *Republic* itself gives no clear answer to it. In a sense, dialectic is the science of being, and the dialectician is 'one who attains a conception of the essence of each thing' by means of the lower sciences.[10] To arrive at this knowledge, Socrates makes 'mathematics, in its fullest development across all its known branches, the basis of the future rulers' training in dialectic'.[11] This point is illustrated by Socrates' analogy in which man begins as a chained prisoner in a cave, where he can only see the shadows of objects presented before him. When he is loosened from his chains and gradually brought out into the open sunlight, he becomes aware of the world as it actually is.[12] Mathematics is, in a sense, Socrates' means of bringing man out of the cave into the contemplation of the sun.[13]

[6] Broadie, Sarah. *Mathematics in Plato's* Republic. Milwaukee: Marquette University Press, 2020, pp. 26-27.
[7] Plato. *Republic*, trans. B. Jowett, in *The Dialogues of Plato*. Oxford: Oxford University Press, 1968, 536d-537a.
[8] Ibid, 537c.
[9] Ibid, 539a.
[10] Ibid, 534b, 533d.
[11] Broadie, p 19.

To see how Socrates can argue thus, it is necessary to understand the value he places on the different mathematical disciplines. In each case, we will see how Socrates is at pains to focus on the subjects' ability to move man from becoming to being.

In discussing the mathematical topics, Glaucon, Socrates's interlocutor, begins with and circles around the utilitarian aspects of the disciplines. Why study arithmetic, geometry, and astronomy? One reason is their use in warfare. Glaucon says that a future ruler should certainly know arithmetic 'if he is to have the smallest understanding of military formations'.[14] Again, Glaucon praises geometry for how it relates to war 'in pitching a camp, or taking up a position, or closing or extending the lines of an army'.[15] Regarding astronomy, Glaucon states how necessary 'the observation of the seasons and of months and years' is to a military general.[16]

But Socrates does not see these materially useful aspects of these fields as the primary reason for their inclusion in the curriculum. At one point, Socrates even pokes fun at Glaucon for his apparent fear of prescribing studies that do not have practical value.[17] So if the purpose of mathematical studies is not in its practical applications, where should one look? Perhaps, to its value as a logically rigorous discipline? Here, again, arithmetic is valuable as a study 'which all arts and sciences and intelligences use in common'.[18] Moreover, talent in calculation tends to make people 'generally quick at every other kind of study'.[19] Arithmetic studies require discipline and tenacity, making it particularly suitable for the curriculum.[20] Regarding geometry, quoting Socrates, 'for the better apprehension of any branch of knowledge, it makes all the difference whether a man has a grasp of geometry or not'.[21]

But even the rigorous method of mathematics is not Socrates's primary concern, for he acknowledges this as one of the 'indirect effects' of the curriculum.[22] According to the English philosopher Myles Burnyeat, Plato excludes from his mathematical disciplines several branches of mathematics which were every bit as rigorous and useful as those branches which he includes. Burnyeat concludes,

> '[t]he great value of mathematics is not practical utility, not transferable skills, not the rigorous procedures of mathematical proof; all these are available from the excluded branches of mathematics.'[23]

---

[12] Plato, 514a et seq.
[13] See Broadie, p. 24 and Burnyeat, M.F., 'Plato on Why Mathematics is Good for the Soul', in T. Smiley, ed., *Mathematics and Necessity: Essays in the History of Philosophy*. Oxford: Oxford University Press, 2001, pp. 21-22.
[14] Ibid, 522e.
[15] Ibid, 526d.
[16] Ibid, 527d.
[17] Ibid.
[18] Ibid, 522c.
[19] Ibid, 526b.
[20] Ibid, 526c.
[21] Ibid, 527c.
[22] Ibid.
[23] Burnyeat, p. 19

That both its utility and rigor should be viewed as indirect effects of mathematical studies is clear because from the outset of the discussion, Socrates has in mind the subjects' ability to 'draw the soul from becoming to being'.[24] In arithmetic, whereas the human eye cannot itself attain essential understanding of the notion of unity, the soul must 'rouse her power of thought…to ask: 'What *is* absolute unity?''.[25] In considering one and all other numbers, the soul is led to contemplate true being.[26] It is in this way that arithmetic achieves its ultimate purpose according to Socrates, for although the man of war or of the marketplace must have knowledge of numbers, it is more useful for the soul herself 'to pass from becoming to truth and being'.[27] In a similar manner, geometry will 'make more easy the vision of the Idea of good' toward which 'all things tend which compel the soul to turn her gaze towards that place where is the full perfection of being'.[28] In geometry, one finds 'knowledge of eternal being, and not of aught which at a particular time comes into being and perishes',[29] which makes geometry essential for inclusion in the curriculum. Even the exquisiteness of the starry heavens is to lead beyond their visible structure towards that which is eternal and unchanging.[30]

Therefore, to Plato, the study of mathematics is, in its fullness, most conducive for contemplation of the truth, leading the soul out of the world of becoming and into the world of being. Through mathematical studies, the future rulers will affect their conversion to the world of being 'in the easiest and quickest manner',[31] for it is in mathematics that one attains knowledge which is not dependent on sense experience and which remain immutably true.

Such a conversion is particularly important for the rulers of the ideal city. To Broadie, Plato sees mathematical studies as the 'cure for…human reason's extreme vulnerability'.[32] The rulers of the ideal city will need to render judgment on matters of ethical and political moment, and Plato needs them to be firmly convinced that in such difficult matters lies the truth. 'By doing [mathematics] they come to know that there are objectively right and wrong solutions to abstract problems and that correct ones can be reached by accurate and imaginatively ingenious reasoning'.[33] In their mathematical studies and debates, the future rulers will be aware that the truth itself is party to their discussion and that they must bear it in mind if they are to arrive at any sound judgment. Hence, the reason Socrates advocates for so much mathematics is not for applications or rigor, per se, but that it is 'an introduction into the whole world of intellect. They learn to believe in that world, in its reality and objectivity, only through total immersion for a long time'.[34]

---

[24] Plato, 521d.
[25] Ibid, 524e.
[26] Ibid, 525a.
[27] Ibid, 525c.
[28] Ibid, 527e.
[29] Ibid, 527b-c.
[30] Ibid, 530b.
[31] Ibid, 518d.
[32] Broadie, p. 41.
[33] Ibid, p. 50.
[34] Ibid, p. 51.

In conclusion, Plato prescribes mathematics primarily for its aid in the development of right reason and for its ability to lead the mind to truth and being. We will now turn to the thought of St. Augustine, where many of the ideas articulated by Plato will be further elucidated and placed into a Christian perspective.

3. Augustine

For Augustine, God is the source and the aim of all knowledge. All of man's knowledge comes ultimately from God, and all truth necessarily leads back to him. In his writings on the nature of man's reason, Augustine places particular importance on the truth of number. Through analyzing number's place in reason's development, we will see the means by which number can lead us to deeper, fundamental truths.

One of Augustine's earliest works following his conversion to Christianity is *De Ordine*, (*On Order*). Augustine's own description of the work is that it is 'an order of studies by which one can proceed from corporeal to incorporeal things'.[35] To that end, Augustine discourses on the liberal arts and their ability to draw the soul towards understanding the things of God. The particular subjects fall into the typical Trivium and Quadrivium topics: Grammar, Logic, Rhetoric, Music, Geometry, Astronomy, and Arithmetic.

Throughout his discourse, Augustine personifies reason as moving through the different disciplines and as noting how in each study, number plays a prime role. In this regard, Augustine goes further than Plato in insisting on the dependence of all the disciplines on numerical ratio. Even within grammar, reason was 'not unmindful of numbers and measure',[36] breaking the words down into their syllables and separating short sounds from long. Within music, 'numeric proportions held sway' and were acknowledged to be immortal.[37] In geometry, reason 'realized that nothing pleased it but beauty; and in beauty, design; and in design, dimensions; and in dimensions, number'.[38] Nothing corporeal pleased the mind so much as what its own intelligence was able to comprehend. Thus were all the liberal disciplines 'presented to reason as numerically proportioned'.[39]

So too, when reason considers the physical world, all that it beholds is number. Number is present in the dimensions and proportion of physical objects, in the motion of the artisan fashioning them, and in the duration of time and physical space in which they are produced.[40] In fact, it is present in the very form of objects themselves, for '[i]f you take a look at changeable reality, you will be unable to grasp it either by the bodily senses or by mental reflection unless it is held together by some numerical determinant, without which it will fall back into nothing'.[41]

---

[35] Augustine. *The Retractions (Retractiones)*. Trans. Mary Inez Bogan, in Roy Joseph Deferrari, ed., *The Fathers of the Church: A New Translation*, vol. 60. Washington, D.C.: The Catholic University of America Press, 1968, 1.3.1.

[36] Augustine. *Divine Providence and the Problem of Evil (De Ordine)*. Trans. Robert Russell, in Ludwig Schopp, ed., *The Fathers of the Church: A New Translation*, vol. 5. New York: Cima Publishing Company, 1948, 2.12.36.

[37] Ibid, 2.14.41.

[38] Ibid, 2.15.42.

[39] Ibid, 2.15.43.

[40] Augustine. *The Free Choice of the Will (De Libero Arbitrio)*. Trans. Robert Russell, in Roy Joseph Deferrari, ed., *The Fathers of the Church: A New Translation*, vol. 59. Washington, D.C.: The Catholic University of America Press, 1968, 2.16.42.

Along with Plato, Augustine acknowledges that the liberal disciplines are studied in part for their practical utility. However, their most beneficial use is for the 'knowledge and contemplation of things' by the use of 'simple and intelligible numbers'.[42] But this is no simple task. To prevent his eyes from being turned aside by corporeal objects, man must seek to unify all he has learned so that he can apply his reason to great matters. In other words, man must use his knowledge of number wisely.

To that end, in a point that he elaborates in his *De Libero Arbitrio (Free Choice of the Will)*, Augustine says that number and wisdom, the proper means by which one contemplates the nature of things, are 'one and the same thing'.[43] Augustine attests to the fact that Scripture relates the two concepts: 'I have gone round—I and my heart—to know and consider, and to search out wisdom and number'.[44] Further, Augustine says that number and wisdom are both part of God's creative wisdom which 'reaches mightily from one end of the earth to the other, and she orders all things well'.[45] While they are, thus, united, number and wisdom play different roles. This point is summarized nicely by Rev. William Most:

> 'Creative Wisdom…gives to all created beings a participation in itself. To rational beings is given wisdom in the full sense of the word. To irrational beings is given number, which is probably to be considered as consubstantial with wisdom, or at least as having the same source.'[46]

Because of the different roles number and wisdom play, Augustine acknowledges that man puts less value on number than on wisdom because the ability to count and number objects is common, whereas the attainment of wisdom is rare. However, since '[w]isdom has endowed all things with number, even the least and those at the lowest confines of the universe,' the truly learned 'behold in the truth itself both number and wisdom and hold both in high esteem'.[47]

In striving to be wise, Augustine is consoled by Scripture, which states that wisdom 'is found by those who seek her. She hastens to make herself known to those who desire her'.[48] In the mathematical ordering of the universe, Augustine sees the imprint of the Wisdom who fashioned it. Just as number is present in all physical objects, we are to see that God is present in giving form to all things. Thus, number, considered in its fullness, leads the soul beyond physical objects to the vision of the artisan at work:

---

[41] *De Libero Arbitrio*, 2.16.44.
[42] *De Ordine,* 2.16.44.
[43] *De Libero Arbitrio*, 2.11.30.
[44] Eccles. 7:25, LXX; quoted in *De Libero Arbitrio*, 2.8.24. Augustine follows the Latin rendering of the Greek Septuagint. The Greek word for 'number' in this verse is 'psēphos', which means pebble. In Augustine's Latin translation, the word is 'numerum'. One could also think of the Latin word '*calculus'* in a similar context.
[45] Wisdom 8:1, RSV; quoted in *De Libero Arbitrio*, 2.11.30.
[46] Most, William G. 'The Scriptural Basis for St. Augustine's Arithmology.' *The Catholic Biblical Quarterly*, 13, No. 3 (July 1951), p. 294.
[47] *De Libero Arbitrio*, 2.11.31.
[48] Wisdom 6:12-13, RSV; quoted in *De Libero Arbitrio*, 2.16.41.

> ‘The artist, too, through the beauty of his work, intimates in a way to the viewer of it that he should not fasten his attention there completely but should so scan the beauty of the artistic work that he will turn his thoughts back fondly upon him who made it.’[49]

Augustine makes the idea of looking beyond number to the artisan at work more concrete in an earlier discussion in *De Libero Arbitrio*. In Book II, Augustine provides a rational proof of God’s existence. To arrive at such a proof, Augustine must surpass that which is greatest in man, namely his ability to reason, ‘for in as much as nothing in man is superior to reason, to transcend reason is to transcend man and consequently to attain God’.[50] Augustine notes that simply finding something that is greater than reason is not itself sufficient to prove that God exists. However, by noting that the senses and reason itself are subject to change,

> ‘If reason sees something eternal and changeless not by any bodily organ…nor by any sense inferior to it, but sees this of itself, and sees at the same time its own inferiority, it will have to acknowledge that this being is its God…or, if there is something higher…that it is God.’[51]

Hence, we see four criteria which Augustine will pursue in the next several chapters of his work. Namely, that he is looking for something that (1) is changeless, (2) is eternal, (3) does not become known through external senses, and finally (4) is greater than reason itself.

In his search for such a thing, it is natural that the discussion turns to truth, and it is especially noteworthy that the first topic that occurs is number. The truth of number is equally available to all who use reason, and the nature of that truth is unchanging and eternal; that seven and three are ten is not a truth today or at this time, but it is always the case.[52] Moreover, knowledge of unity in number is not something that comes through the senses, for in as much as every body admits of parts (a left and a right, a top and a bottom, etc.), ‘we acknowledge that no bodily reality is one, truly and simply’.[53] But even in the acknowledgment that no corporeal thing is one comes knowledge of oneness itself. Moreover, in asserting a law or theorem of number, the senses cannot perceive that it holds for all numbers, for they are innumerable. Hence, it is that one apprehends number by an ‘inner light of which the senses have no knowledge’.[54]

Thus, number is an example of a changeless, eternal truth which is not apprehended through the senses. Augustine’s argument establishes that ‘there exists unchangeable truth that embraces all things that are immutably true’ whether they be truths of number or whether they be truths of a different sort.[55] This class, which could be called truth itself, is what Augustine says is superior to our minds and to reason. His argument is one by trichotomy. Truth cannot be less than reason, for we make judgments according to the truth, but we do not make judgments about the truth. Augustine likens it to the mind making discoveries in the light of the truth rather than as an examiner making corrections to the truth.[56] Moreover, it cannot be equal to reason, for then the truth would be changeable just as our minds grow or diminish in capacity. Therefore, since it is neither less than nor equal to our minds, the truth must be greater than our minds.

---

[49] *De Libero Arbitrio*, 2.16.43.
[50] Gilson, Etienne. *The Christian Philosophy of Saint Augustine*. Providence: Cluny Media, 2020, p. 19.

[51] *De Libero Arbitrio*, 2.6.14.
[52] Ibid, 2.8.20-21.
[53] Ibid, 2.8.22.
[54] Ibid, 2.8.23.
[55] Ibid, 2.12.33.

And there it is, the goal to which Augustine has been working: ‘For if there is anything more excellent [than this truth], then this is God; if not, then truth itself is God. In either case, you cannot deny that God exists’.[57] In this truth should man find his happiness and, in time, be brought into the fullness of wisdom. So far has reason alone brought Augustine, but he sees further that in our search, we hear Truth Himself speaking to us: ‘If you continue in my word, you are truly my disciples, and you will know the truth, and the truth will make you free’.[58] Therefore, to Augustine, the truth of numbers leads us to see that there is truth unshakeable. We see through it that God exists, the craftsman who fashioned things to be such.

4. Conclusion

In his address, Shafarevitch laments the loss of a global aim of humanity’s cultural activity as indicated in mathematics’ own predicament. Such a state can quickly lead to a world which puts too much emphasis on progress and reason to the detriment of humanity. In his book *Truth and Tolerance*, Joseph Ratzinger (the future Pope Benedict XVI) says that while amassing ever more knowledge and precision, reason left to itself ‘no longer offers any perspective on the fundamental questions of mankind’.[59] The only solution to the problem is one of unification in aim: ‘[R]eason and religion will have to come together again…neither of them, it is clear, can be saved unless God reappears in a convincing fashion’.[60] This is certainly a view with which Shafarevitch would heartily agree.

In his discourses, Augustine follows his classic formula of believing so as to understand, and we see in his approach an integration of faith with reason. It is in this way that Shafarevitch’s conception of the fundamental religious meaning of mathematics can be realized. Speaking of the Pythagoreans, Shafarevitch observes that ‘[b]y revealing the harmony of the world as expressed in the harmony of numbers it provided a path leading towards a union with the divine’.[61] Finding inspiration from the Pythagoreans, Plato saw the ability of mathematics to raise the soul out of the world of becoming into the world of being. Augustine picked up this thread and placed it into a Christian perspective. To be clear on Augustine’s perspective, number is not to be the end of the story, for

> ‘whoever esteems [number] so highly that he wants to speak boastfully of himself among unlearned men, instead of trying to learn the source of the truth of things which he has seen are true…may seem to be erudite, but he can by no means be considered wise.’[62]

---

[56] Ibid, 2.12.34.
[57] Ibid, 2.15.39.
[58] John 8:31-32 RSV; quoted in *De Libero Arbitrio*, 2.13.37.
[59] Ratzinger, Joseph. *Truth and Tolerance*. San Francisco: Ignatius Press, 2004, p. 143.
[60] Ibid, p. 144.
[61] Shafarevitch, p. 9.
[62] Augustine. *Christian Instruction (De Doctrina Christiana).* Trans. John Gavigan, in Ludwig Schopp, ed., *The Fathers of the Church: A New Translation*, vol. 4. New York: Cima Publishing Company, 1947, 2.38.57.

Indeed, man's knowledge must be unified. For Augustine the philosopher, 'both in analyzing and in synthesizing, it is oneness that I seek, it is oneness that I love'.[63]

Oneness and unity form the heart of Augustine's vision of the purpose and use of number. Number is present in all material objects, but to Augustine, 'though it has parts, a body imitates the Unity that is God, who is Being Itself'.[64] And in this imitation, we see Wisdom hastening to make itself known and to draw us into union with itself. As Augustine says, 'What about every kind of love? Does it not wish to become one with what it is loving? And, if it reaches its object, does it not become one with it?'.[65]

To conclude, Augustine's treatment of number also carries with it a moral component. To the soul who has encountered the truth and nature of numbers and their abode in eternal truth, 'it will appear manifestly unfitting and most deplorable that it should write a rhythmic line and play the harp by virtue of this knowledge and that its life and very self…should nevertheless follow a crooked path and…be out of tune by the clangor of shameful vices'.[66] So by ordering our lives appropriately and applying ourselves diligently to our studies, with God's help, we may have power to comprehend 'what is the breadth and length and height and depth, and to know the love of Christ which surpasses knowledge'.[67]

---

[63] *De Ordine,* 2.18.48.
[64] Most, p. 294.
[65] *De Ordine,* 2.18.48.
[66] *De Ordine,* 2.19.50.
[67] Eph. 3:18-19, RSV.